\documentclass[oneside,english]{amsart}
\usepackage[T1]{fontenc}
\usepackage[latin9]{inputenc}
\usepackage{url}
\usepackage{amsthm}
\usepackage{amssymb}
\usepackage{todonotes}
\usepackage{hyperref}

\makeatletter
\numberwithin{equation}{section}
\numberwithin{figure}{section}
\theoremstyle{plain}
\newtheorem{thm}{\protect\theoremname}
\theoremstyle{definition}
\newtheorem{defn}[thm]{\protect\definitionname}
\theoremstyle{definition}
\newtheorem{problem}[thm]{\protect\problemname}
\theoremstyle{remark}

\theoremstyle{plain}
\newtheorem{question}[thm]{\protect\questionname}
\theoremstyle{definition}

\newtheorem{conjecture}[thm]{Conjecture}
\include{mydefs}

\newcommand{\Ao}{\mathbb{A}^1}

\newcommand\GW{\mathrm{GW}}
\DeclareMathOperator\Poly{Poly}
\DeclareMathOperator\Hol{Hol}
\DeclareMathOperator\Rat{Rat}
\DeclareMathOperator\Mor{Mor}
\DeclareMathOperator\Mod{Mod}

\usepackage{mathtools}
\DeclareMathOperator\GL{GL}
\DeclareMathOperator\PGL{PGL}
\DeclareMathOperator\Hom{Hom}
\DeclareMathOperator\Char{char}
\DeclareMathOperator\Sel{Sel}
\newcommand{\PP}{\mathbb{P}}
\newcommand{\CP}{\mathbb{CP}}
\newcommand{\CC}{\mathbb{C}}
\newcommand{\RR}{\mathbb{R}}
\newcommand{\FF}{\mathbb{F}}
\newcommand{\QQ}{\mathbb{Q}}
\newcommand{\ZZ}{\mathbb{Z}}
\renewcommand{\Mc}{\mathcal{M}}

\makeatother

\usepackage{babel}
\providecommand{\definitionname}{Definition}
\providecommand{\examplename}{Example}
\providecommand{\problemname}{Problem}
\providecommand{\questionname}{Question}
\providecommand{\remarkname}{Remark}
\providecommand{\theoremname}{Theorem}

\title{Problems in Arithmetic Topology}
\author{Claudio G\'omez-Gonz\'ales}
\author{Jesse Wolfson}

\address{ Department of Mathematics, University of California-Irvine}
\email{claudijg@uci.edu}
\email{wolfson@uci.edu}

\thanks{JW is supported in part by NSF Grant Nos. DMS-1811846 and DMS-1944862. Support for the workshop was provided by PIMS and by NSF Grant No. DMS-1856737. }

\begin{document}

\maketitle

Three problem sessions were hosted during the workshop in which participants proposed open questions to the audience and engaged in shared discussions from their own perspectives as working mathematicians across various fields of study. Participants were explicitly asked to provide problems of various levels of difficulty, with the goal of capturing a cross-section of exciting challenges in the field that could help guide future activity.  The problems, together with references and brief discussions when appropriate, are collected below into three categories: 1) topological analogues of arithmetic phenomena, 2) point counts, stability phenomena and the Grothendieck ring, and 3) tools, methods and examples.

\subsection*{Acknowledgements}
We thank all of the people who posed problems, both for their original suggestions and for their help clarifying and refining our write-ups of them. We thank Christin Bibby, Aaron Landesman, Lior Silberman, and Bena Tshishiku for transcribing the problem sessions in real time. We thank Benson Farb for helpful commnts. We thank all the participants who commented and asked questions during the problem sessions. Last, we thank Alejandro Adem, Craig Westerland and Melanie Wood for co-organizing the workshop and we thank the PIMS staff for their hospitality and support in organizing the workshop.

\section{Topological Analogues}
%1
\subsection{Craig Westerland}

This problem concerns the homology of random topological objects and the Cohen--Lenstra distribution, originally introduced to heuristically understand experimental observations about the class groups of number fields \cite{CL}.  The Cohen--Lenstra distribution on finite (abelian) groups posits that the probability of each group $G$ being in the support of the distribution is inversely proportional to $|G|^n |\operatorname{Aut}(G)|^m$ for some $n, m \geq 0$. Since its first use, it has appeared throughout number theory and in various models of random integral matrix cokernels; see, for example, \cite{CKLPW, EVW, koplewitz, wood15}.

Working one prime $p$ at a time, one can compute the $p$-power torsion and torsion-free part of the integral homology of a given chain complex by tensoring with the $p$-adic integers $\ZZ_p$ before taking homology. If one further specifies the total homological rank $r$ of such complexes, the Cohen--Lenstra distribution specifies a discrete measure on the resulting module $M = \ZZ_p^r \times N$ (for $N$ an abelian $p$-group). We can then ask the following question:

\begin{conjecture} Let $X_\bullet$ be a random simplicial complex. Does the torsion $H_\ast(X_\bullet)_{tors}$ have a Cohen--Lenstra distribution? 
\end{conjecture} 

This question was originally posed by Matt Kahle, Frank Lutz, Andrew Newman and Kyle Parsons, along with experimental evidence supporting the conjecture, in \cite{KLNP}. Westerland raised this standing problem and also discussed an analogue for random manifolds where one can introduce a torsion linking pairing
\[ \begin{split} H^i(X;\ZZ)_{tors} \times H^{\dim X - i - 1}&(X;\ZZ)_{tors} \stackrel{\cup}{\longrightarrow} H^{\dim X - 1}(X;\ZZ)_{tors} \\
	& \stackrel{\delta}{\longrightarrow} H^{\dim X}(X;\QQ/\ZZ) \cong \QQ/\ZZ.
\end{split} \]
Delaunay \cite{delaunay} introduced heuristics in the context of a bilinear pairing \[\mu: G \times G \to \QQ/\ZZ,\] originally in the symplectic context of Tate--Shafarevitch groups but later expanded by others (e.g. \cite{BKLPR, wood14}) to additional phenomena, where $\operatorname{Aut}(G)$ is replaced with $\operatorname{Aut}^\mu(G)$ of automorphisms preserving the pairing.

\begin{question} Are there Delauney heuristics in the random manifold case?
\end{question}

Andy Putman mentioned useful models for random $3$-manifolds (e.g., Dunfield--Thurston manifolds via random classes in the mapping class group $\Mod_g$) that might help approach this problem.

%2
\subsection{Daniel Litt}
Motivated by the Frobenius action on the pro-$p$ Galois group of a curve over $\FF_p$, consider a random automorphism $\varphi$ of a free pro-$p$ group $F$, and construct a random Weil group 
\[ W_\varphi \coloneqq F \rtimes \langle \varphi\rangle.\]

What is the representation theory of $W_\varphi$? Do the predictions from geometric Langlands theory hold? For example:
\begin{problem}  
	Given an irreducible continuous representation $\rho: W_\varphi \to \GL_n\left(\FF_p((t))\right)$, does the restriction $\rho|_F$ always have finite image?  
\end{problem}
For $F$ the geometric \'etale fundamental group of a smooth projective curve $X$ over $\FF_q$, the analogous problem was posed by de Jong \cite{dJ} and proven by Gaitsgory \cite{Ga}. Similarly, one could ask:
\begin{problem}
	Let $L/\QQ_p$ be a finite extension, and let $\rho\colon W_\varphi\to \GL_n(L)$ be a continuous irreducible representation.  For almost every $\varphi$ (i.e. with probability 1), does $H^1_{cts}(F;\rho\otimes\rho^\vee)^\varphi=0$?
\end{problem}
Viewing $H^1_{cts}(F;\rho\otimes\rho^\vee)^\varphi$ as the tangent space to $\rho$ in the character variety of $F$, we can understand the previous problem as a ``cohomological rigidity'' statement.  Like the first problem, this problem is a theorem for $F$ the geometric \'etale fundamental group of a smooth projective curve $X$ over $\FF_q$: indeed, a theorem of L. Lafforgue \cite{Lafforgue} implies that $\rho$ is pure of some weight, so $\rho\otimes\rho^\vee$ is pure of weight 0.  Deligne's bounds \cite{Deligne} therefore imply that $H^1_{cts}(F;\rho\otimes\rho^\vee)$ is pure of weight 1, and hence the invariants of Frobenius are trivial.

More generally, both of the above problems arise from trying to understand in what sense the action of Frobenius on the geometric \'etale fundamental group of a smooth projective curve over $\FF_q$ is ``special''?  Affirmative answers to both of the questions above would suggest that it isn't!

\section{Point counting, Stability Phenomena and the Grothendieck Ring}

%3
\subsection{Akshay Venkatesh} In many interesting number theory problems (e.g. counting number fields) one has not only a main term in the asymptotic count, but a secondary term or more.  For example, the number of cubic fields of discriminant up to $X$ is
\[
a X + b X^{5/6} +\text{lower order terms}
\] 
We have very little understanding of these lower order terms, and they are not just of theoretical interest: when one tries to verify the conjectures numerically, one finds that the secondary terms are dominant in the computational range. 

So the question, following the framework of the Weil conjectures, is:
\begin{problem}
	What is the topological meaning of secondary terms appearing in assymptotic counts in number theory?
\end{problem}
They do not correspond to stable cohomology classes (these are the main terms), but to some kind of slightly weaker structure, which is still much better behaved than cohomology near middle degree. 

%4
\subsection{Zachary Himes}

The recent work of Galatius--Kupers--Randal-Williams \cite{GKRW1, GKRW2} develops and applies a homology theory based on $E_k$-algebras, originally introduced to study $k$-fold loop spaces \cite{may}, with the goal of showing results ``beyond homological stability'' which they call \emph{secondary stability}. For example, in \cite{GKRW2}, $E_k$-homology is used to measure the failure of homological stability itself as a stable phenomenon in the context of mapping class groups, and in \cite{RW}, this perspective is used to explain the homological stability phenomena which underlie the breakthrough work of Ellenberg--Venkatesh--Westerland \cite{EVW}. 

\begin{question} Building off \cite{EVW}, can one prove Hurwitz spaces have stability (for more general groups than the dihedral groups)?  Secondary stability, or higher-order-stability? What can one say  about the homology of Hurwitz spaces after restricting to subsets of connected components? 
\end{question}

%5
\subsection{Daniel Litt}

\begin{problem} Fix $n, d,$ and $q$. How many supersingular hypersurfaces $X \subset \PP^n$ exist of degree $d$ over $\FF_q$? In particular, does such a hypersurface exist for every degree, dimension, and characteristic? Moreover, the set of such varieties are the $\FF_q$-points of some supersingular locus $V$---what is the class $[V]$ in the Grothendieck ring? 
\end{problem}

A great deal is known about the supersingular locus in particular cases. See, for example, the work of Oort and Li \cite{lioort} on moduli spaces of supersingular abelian varieties, especially on the case of principally polarized abelian varieties $S_{g,1}$.

%6
\subsection{Aaron Landesman}

Cohen--Lenstra heuristics can be used to describe distributions of class groups of certain global fields. In many ways, the Selmer groups of elliptic curves behave similarly to class groups. As such, Landesman provided two conjectures---a number-theoretic version and a topological version---relating topological techniques to the average size of Selmer groups of elliptic curves over function fields. We describe the latter here, framed as a question of homological stability, which would likely imply the number-theoretic conjecture. For more details, see \cite{landesman}.

\begin{defn}
	The \emph{moduli space of Weierstrass models} is given by
	\[ \begin{split}
		\mathcal{W}_d \coloneqq \{ \left( a(S,T),b(S,T) \right) \in \CC[S,T]^2 \mid & \text{ homogeneous degree $4d$} \\
		& \quad \text{and $6d$, respectively}\} \subseteq \CC^{10d+2}. 
	\end{split} \]
	A point $W = (a(S,T),b(S,T)) \in \mathcal{W}_d$ gives rise to an elliptic surface $S_W$ via
	\[ S_W \coloneqq \{ ([S:T],[X:Y:Z]) \in \CP^1 \times \CP^2 \mid Y^2Z = X^3 + a(S,T)XZ^2 + b(S,T)Z^3 \}, \]
	together with the natural projection map $\varphi_W: S_W \to \CP^1$. The \emph{universal family of Weierstrass models}, $\mathcal{UW}_d$, is the space parameterizing models $W$ and a point $p$ on the corresponding elliptic surface $S_W$. That is,
	\[ \mathcal{UW}_d \coloneqq \{ (W,p) \in \mathcal{W}_d \times (\CP^1 \times \CP^2) \mid p \in S_W \}. \]
	There are natural projection maps
	\[ \mathcal{UW}_d \stackrel{f}{\longrightarrow} \CP^1 \times \mathcal{W}_d \stackrel{ g}{\longrightarrow} \mathcal{W}_d. \]
	
	The Selmer space is constructed as $\Sel_n^d(\CC) \coloneqq R^1 g_*( R^1 f_*(\mu_n))$.
\end{defn}

\begin{question}
	For fixed $n$, do $\Sel_n^d(\CC)$ satisfy homological stability as $d \to \infty$?
\end{question}

In particular, Landesman provided a more precise prediction:

\begin{conjecture}
	There are constants $A$ and $B$ depending on $n$ so that
	\[ \dim H_i(\Sel_n^d(\CC); \QQ) = \dim H_i(\Sel_n^{d+1}(\CC);\QQ) \]
	whenever $d \geq Ai+B$.
\end{conjecture}

%7
\subsection{Will Sawin}
Continuing the paradigm put forward by Vakil--Wood \cite{VW}, Sawin proposed further comparison of three types of convergence for families of varieties:
\begin{enumerate}
	\item Homological stability, i.e. isomorphisms between low-degree singular cohomology groups (equivalently, high-degree compactly-supported cohomology groups), ideally in some way respecting arithmetic structures such as Galois representations, Hodge structures, etc.
	\item Point counts, in particular after normalizing by a factor of $q^{-\dim}$.
	\item Convergence in the Grothendieck ring under the dimension filtration; see, for example, \cite{VW}.
\end{enumerate}

In particular, Sawin suggested fixing a degree $d$ and considering a sequence $H_n$ of smooth degree $d$ hypersurfaces in $\PP^n$ as follows:

\begin{question}
	Does $\displaystyle \lim_{n \to \infty} \frac{[H_n] - [\PP^{n-1}]}{[\mathbb{A}^{n-1}]}$ exist? Is the limit zero?
\end{question}

This convergence makes sense for many objects that came be obtained via maps out of the Grothendieck group (i.e., Hodge structures), so a negative result would show the Grothendieck ring is richer than these. In the setting of homological stability, this question reduces to the Lefschetz hyperplane theorem; with respect to counting points over $\FF_q$, the limit can be controlled when $\sqrt{q}$ is larger than $d-1$ using Lefschetz together with the Weil conjectures. While any progress towards answering this question is interesting in its own right, there is also an application to proving Grothendieck ring analogues of the Browning--Sawin circle method \cite{browningsawin}.

%8
\subsection{Ronno Das}

In 1849, Cayley and Salmon \cite{cayley} showed that every smooth cubic surface over an algebraically closed field contains exactly 27 lines. If we instead let $X$ be a smooth cubic surface over $k = \RR$ or $\FF_q$, for example, the number of lines contained in $X$ and defined over $k$ can be strictly less than 27. In fact, as shown by Das \cite{das}, the average number of lines on a cubic surface over $\FF_q$ is exactly 1 for all but finitely many characteristics. 

More generally, as suggested in correspondence with Ravi Vakil, one can consider lines on del Pezzo surfaces of degree $1 \leq d \leq 9$, where $d=3$ is the case of cubic surfaces. For more details see, for example, \cite{BG} and \cite{josef}.

\begin{problem} 
What is the average number of $\FF_q$-lines of a degree 1 del Pezzo surface over $\FF_q$? % apparently Bergvall has a preprint for d=2, but I don't know if I should ask him about it? % haven't heard back from Bergvall
\end{problem}

The average number of lines on a del Pezzo surface defined over $\FF_q$ is $1$ for degree $d \geq 3$, except when $d=7$ where the average is 2. Moreover, the $d=2$ case is closely related to Bergvall's computations \cite{bergvall}.

\begin{problem} 
Is there a uniform proof for all degrees $1 \leq d \leq 9$?
\end{problem}

%9
\subsection{Jesse Kass}

When $q$ is odd or when $k = \RR$, Segre \cite{segre} observed that a line $\ell$ on a cubic surface $X$ comes equipped with a distinguished involution $L$ that is necessarily \emph{hyperbolic} or \emph{elliptic} as an element of $\PGL_2$, determined by whether or not the fixed points of $L$ are defined over $k$. 

\begin{defn} A line $\ell$ on a cubic surface $X$, both defined either over $\RR$ or $\FF_q$ for $q$ odd, is said to be \emph{hyperbolic} (respectively, \emph{elliptic}) if the distinguished involution $L$ (see \cite{segre}) is hyperbolic (respectively, elliptic). \end{defn} 

See the papers of Finashin--Kharlamov \cite{finashinkharlamov}, Okonek--Teleman \cite{okonekteleman}, and Kass--Wickelgren \cite{KW17} on how this distinction gives rise to restrictions for the number of lines on cubic surfaces defined over $\RR$ and $\FF_q$ with odd characteristic.

\begin{problem} Let $X$ be a cubic surface over $\FF_q$, a finite field of odd characteristic, and with 27 lines defined over $\FF_q$. Let $H$ and $E$ be the number of hyperbolic and elliptic lines, respectively.
\begin{enumerate}
\item What pairs $(H,E)$ can appear?
\item Statistically, how often do each appear?
\item What subgraphs of the intersection graph of lines can be given by hyperbolic and elliptic lines?
\end{enumerate}
\end{problem}

%10
\subsection{Isabel Vogt}

\begin{problem} Let $X$ be a smooth cubic hypersurface defined over $\FF_q$ and fix $n \leq q+1$ points from $X$. Does there exist an algebraic map $f: \PP^1_{\FF_q} \to X$ such that the $n$ points are contained in the image? If so, what degree is necessary?
\end{problem}

As proved by Koll\'ar \cite{kollar}, the answer to the first question is yes when $q \geq 8$, but the degree of $f$ grows rapidly. In general, for $X$ a fixed hypersurface with $g: Z \to X$ a map on some finite subscheme $Z$ of $\PP^1$, one can consider the space
\[ \Mor_d(\PP^1,X;g) \coloneqq \left\{ f \in \Mor_d(\PP^1,X) : f|_Z = g \right\} \]
of maps extending $g$ to all of $\PP^1$. The preceding question becomes that of the $\FF_q$-points of $\Mor_d(\PP^1,X;g)$. In this spirit, Vogt offered the following general problem:

\begin{problem}\label{extend_top} Asking these questions over $\CC$, how does the topology of the associated spaces vary as $d$ grows?
\end{problem}

As a related example, fixing a degree $d$ holomorphic function $g: \CP^{m-1} \to \CP^n$, Mostovoy \cite{Mostovoy03, Mostovoy12} considered the space
\[ \Mor_d(\CP^m,\CP^n;g) \coloneqq \left\{ f: \CP^m \to \CP^n | f \circ i = g \right\} \]
for some fixed embedding $i: \CP^{m-1} \hookrightarrow \CP^m$, and exhibited a map
\[ \Mor_d(\CP^m,\CP^n;g) \to \Omega^{2m} \CP^n \]
that induces isomorphisms in homology up to degree $d(2n-2m+1)-1$. Since the latter space has a homotopy type naturally invariant of $d$, this result (and subsequent generalizations into toric varieties \cite{MMV}) gives a homological stability result in the spirit of problem \ref{extend_top}.

%11
\subsection{Vlad Matei}

\begin{problem} Consider the variety of monic square-free polynomials,
\[ \Poly_n(k) \coloneqq \{ (a_{n-1},\dots,a_0) \in k^n \mid t^n + a_{n-1} t^{n-1} + \cdots + a_0 \text{ is square-free} \}. \]
%When $k = \CC$, this variety is the classifying space of the Artin braid group. 
Fixing $n$ and $k = \FF_q$, what is the distribution of the $a_i \in \FF_q$ for $i \leq n-2$?
\end{problem}

In the $q \gg 1$ regime, the coefficients are well-known to be uniformly distributed.
On the other hand, fixing $q$, one can see immediately that $a_{n-1}$ is uniformly distributed simply by elementary substitutions. The question of fixing some of the coefficients is of particular interest to computer scientists; see, for example, \cite{granger} for a history of the problem. Matei outlined two distinct ways to study the aforementioned problem:

The first is by dualizing---considering the discriminant hypersurface $\Delta \subset \mathbb{A}^n$ and picking out $\FF_q$-points by their coordinates. These, in turn, can be thought of as hyperplane sections of $\Delta$. The downside of this approach is that we would need an explicit formula for $\Delta$ in terms of the coefficients, which becomes computationally infeasible even in small degrees.

The second approach is by considering hypersurfaces determined by the symmetric sums of the roots $\{z_1,\dots,z_n\}$. Writing $\sigma_k$ for the $k$th symmetric sum, one could consider hypersurfaces of the form $\sigma_k = a \in \FF_q$ on the open locus of distinct roots. In this case, the objective is to understand the $S_n$-action on these hypersurfaces and the associated cohomology groups, in order to count $\FF_q$ points on the quotient.

%12
\subsection{Oishee Banerjee}

While many of the results concerning the cohomology of moduli spaces of polynomials concentrate on data encoded in their zeroes, i.e. configuration spaces on $\CC$, one can also proceed by studying ramification. Such work is in the spirit of long-standing open problems concerning the topology of Hurwitz spaces. 

Indeed, the irreducibility of the Hurwitz space is a classical result proved in \cite{clebsch}, where the topology of its subvarieties corresponding to specific ramification loci is almost completely unknown. Banerjee \cite{banerjee} studied the stable cohomology of these Hurwitz spaces satisfying certain conditions. In addition her work shows that the \'etale cohomology does not stabilize in positive characteristic, which is in contrast to comparable stability results (see \cite{FW, EVW}.) 

\begin{problem} Fix a prime power $q$ and integer $n$ with $\Char(\FF_q) > n+1$. For which $f(t) \in \Poly_n(\FF_q)$ does the anti-derivative 
\[ F(t) \coloneqq \int_0^t f(x) dx \]
give a simply-branched mapping $F: \Ao \to \Ao$? Alternatively, how many simply-branched maps $F: \Ao \to \Ao$ over $\FF_q$ of degree $n+1$ have $\tfrac{dF}{dt} \in \Poly_n(\FF_q)$?
\end{problem}

%13
\subsection{Patricia Hersh}

This question is motivated by the task of moving from studying the configuration space of $n$ distinct points on a manifold, either entirely labeled or unlabeled, to instead looking at situations where the $n$ points are partitioned into groups of points that cannot be distinguished from each other but can be distinguished from the other groups of points.

If $A$ is a hyperplane arrangement one can count the number of $\FF_q$ points in the complement $X$ of $A$ using the M\"obius function $\mu_A$ of the intersection poset $L_A$ (see, for example, \cite{geometriccombinatorics}). In short: \begin{equation}\label{complement count} \# X(\FF_q) = \sum_{y \in L_A} \mu(\hat{0},y) q^{\dim A - \operatorname{rk} y}. \end{equation}
Hersh--Kleinberg introduced a multiplicative deformation of the M\"obius function, $\mu'$, defined somewhat similarly to $\mu_A$ but (1) replacing the recursion for $\mu_A$ by a recursion that forces $\mu'$ to be multiplicative, and (2) replacing the intersection poset $L_A$ in the counting formula above by the multiset partition poset ordered by refinement. In particular, $\mu'$ records multiplicities of incidences and thus behaves more elegantly in various contexts (see \cite{HK} for details).

\begin{question} Replacing $\mu$ with $\mu'$ in (\ref{complement count}) and replacing $L_A$ by the multiset partition poset, does the resulting formula have arithmetic meaning?
\end{question}

Additional references are \cite{garton} and \cite{hersh}, the former remarked by Melanie Wood.

%14
\subsection{Joseph Gunther}

\begin{question}
	Let $\mathcal{C}(a,b)$ be moduli space of smooth complex curves in $\PP^1 \times \PP^1$ of bi-degree $(a,b)$. Is there a homological stability for $\mathcal{C}(a,b)$ as $b \to \infty$ with $a$ fixed?
\end{question}

This problem is solved over $\FF_q$ by \cite{EW} and a motivic treatment is found in \cite{BH}.

%15
\subsection{Vlad Matei}

In 1984, Kani--Rosen \cite{KR} studied relations between idempotents in the algebra
of rational endomorphisms of a fixed abelian variety. In particular, if a group $G$ covered by subgroups $H_1,\dots,H_r$ acts on a smooth variety $X$ over $\QQ$, we obtain idempotents $\epsilon_{H_i}$ given by
\[ \epsilon_{H_i} \coloneqq \frac{1}{H_i} \sum_{h \in H_i} h \]
in the group algebra $\QQ[G]$ that come with linear relations of the form
\[ \sum_{i=1}^n a_i \epsilon_{H_i} = 0. \]

\begin{question}
	Do these transfer to relations in the Grothendieck ring? That is, are there $a_i \in \ZZ$ such that
	\[ \sum_{i=1}^r a_i [X/H_i] = 0? \]
\end{question}

Matei observed that this holds for the action of $G = S_3$ on affine space $\mathbb{A}^3$, but Daniel Litt commented that he expects this to be false mod $p$ generically.

\section{Developing tools, methods, and new examples}
%16
\subsection{Will Sawin}

Sawin's talk outlined joint work with Tim Browning on a geometric version of the Hardy--Littlewood circle method, used to compute the compactly supported cohomology of the space of rational curves on a smooth hypersurface (see \cite{browningsawin}). 

\begin{question}
	Can we apply these methods to other cases, such as to maps
	\begin{enumerate}
		\item from higher genus curves?
		\item into projective varieties?
		\item into complete intersections?
		\item with local conditions?
		\item with Hodge structure modules? 
	\end{enumerate}
\end{question}

%17
\subsection{Michael Fried} 
For $n\ge 4$, two {\em moduli} spaces associated to $\PP^1_z$ appear often:
\begin{center}
	The moduli of ordered (distinct) points $\Mc_{0,n}:=((\PP^1_z)^n\setminus \Delta_n)/\PGL_2(\CC)$, and\\
	The moduli of unordered (distinct) points $J_n:=(\PP^n\setminus D_n)/\PGL_2(\CC)$
\end{center}
where $\Delta_n$ is the ``fat diagonal'', i.e. the locus of points with repetitions, and $D_n$ is the discriminant (treating the homogeneous coordinates in $\PP^n$ as the coefficients of a degree $n$ polynomial in one variable). For $n=4$, these are the (open) classical lines $\PP^1_\lambda\setminus\{0,1,\infty\}$ and $\PP^1_j\setminus\{\infty\}$. The fundamental group of $\Mc_{0,n}$ mod outer automorphisms, $\pi_n$, is the pure (spherical) braid group on $n$ strands mod its center $\ZZ/2\ZZ$. Note that $\pi_4$ is a free group on two generators, $F_2$.

The Grothendieck--Teichm\"uller group is a conjectured description (cf. \cite{GrEsquisse, IharaICM, DrGT}) of the absolute Galois group $G_{\QQ}$ given by its action on the profinite fundamental group(oid)s of $\Mc_{0,4}$ and $\Mc_{0,5}$. Even if correct, Grothendieck's conjecture is not useful without a way to name elements in the action of $G_{\QQ}$ (cf. \cite[Question 2, p. 6]{IharaICM}).

Ihara replaced $\pi_n$ by $\pi_n^{nil}$, its nilpotent completion, as a way to name the elements in their action; note that this won't capture all of $G_{\QQ}$. He, along with Anderson, documented aspects  of this for explicit collections of covers (e.g. \cite{IharaAnderson1,IharaAnderson2}).  Foremost was his identification of the second commutator quotient of $G_{\QQ}$ (the 2-step nilpotent quotient), using the 3, 4, and 5-cycle relations developed by Drinfeld (and apparently known to Grothendieck). This generalized the Kronecker-Weber description fo the abelian quotient of $G_{\QQ}$ as the (Galois group over $\QQ$ of the) cyclotomic closure of $\QQ$ (see \cite{IharaKW}). 

\begin{problem}
	Go beyond the nilpotent quotient $\pi_n^{nil}$ to give a naming scheme for elements in $G_{\QQ}$.
\end{problem}
 
 We propose that a moduli-theoretic approach is possible using reduced Hurwitz space covers (see \cite{Fried} and \cite{BF,FV}). These covers of $J_n$ -- defined by finite groups and conjugacy classes withing these group, and parametrizing branched covers of $\PP^1$ with fixed branching data and monodromy group -- pull back to unramified covers of $\Mc_{0,n}$ whose components and their fields of definition are controlled by an explicit braid group action.  For example, for $n=4$, components, genuses and cusps of these covers of $\Mc_{0,4}=\PP^1_{\lambda}\setminus\{0,1,\infty\}$ are efficiently computable. We propose that these Hurwitz moduli spaces give a way to identify the $G_{\QQ}$ action.  Using versions of Deligne's tangential base points \cite{DeligneDroite} based in $J_n$ rather than $\Mc_{0,n}$ (significant even for $n=4$)  thereby gives tests relating Drinfeld's relations to the fields of definition of components of Hurwitz sapce. 
 
{\bf Special case:} even the case when the Hurwitz space pullbacks to $\Mc_{0,n}$ are nilpotent covers should be doable and significant (e.g. in shedding new light on Ihara's results).

%18
\subsection{Jesse Kass}

The following was discussed in the context of developing an enriched discriminant or conductor using $\Ao$-homotopy theory, in the context of the mini-course by Kirsten Wickelgren and work by Marc Levine; see, for example, \cite{levine17}.

Consider a generically separable non-constant Galois cover $f: Y \to X$ of smooth projective curves with Galois group $G$. Associated to $df$ there is a certain Euler class, the \emph{enriched different}, which is related to the ``enriched Riemann--Hurwitz formula'' due to Levine. In classical theory, the different is related to the discriminant and the conductor; see, for example, \cite{neukirch} for an overview.

\begin{question}
	Is there an \emph{enriched conductor} associated to $y \in Y$ and $x = f(y)$? 
\end{question}

\begin{defn}
	The classical conductor is defined via the Artin character (see \cite{greenbergserre}). In this context, where $G_y \coloneqq \operatorname{Stab}_G(y)$, we have an \emph{enriched Artin character}
	\[ \begin{gathered} 
		G_y \to \GW(k) \\
		\sigma \mapsto \text{local Lefschetz trace of } \sigma,
	\end{gathered} \]
	where $\GW(k)$ is the Grothendieck-Witt group associated to the field $k$.
\end{defn}

\begin{question}
	Is the enriched Artin character the character of some map?
\end{question}

More generally, the task here is to continue enriching Serre's book \cite{greenbergserre}. 

\subsection{Daniel Litt}

Following the work of Marc Levine, Jesse Kass, Kirsten Wickelgren, and others (see, for example, \cite{KW17,KW19,levine18}), Daniel Litt asked about creating enriched Gromov-Witten invariants.

\begin{question} 
	Can we develop enriched Gromov-Witten invariants in the non-enumerative case? Can we relax the conditions, e.g. for non-smooth varieties? 
\end{question} 

Kass remarked that there are computations over the real numbers, for example \cite{GZ,ST,welschinger}, which might give us some hint on how to define them.

%19
\subsection{Hannah Knight}

In 1979, based on intuition from Morse theory, Segal \cite{segal} proved that the map
\[ \Rat_d^n(\CC) \hookrightarrow \Omega^2 \CC\PP^n \]
is a homotopy equivalence through dimension $(2n-1)d$, where $\Omega^2 \CC\PP^n$ is the space of based continuous maps $\CC\PP^1 \to \CC\PP^n$ with the compact-open topology. %and
%\[ \Rat_d^n(\CC) \coloneqq \left\{  (f_0(t),\dots,f_n(t)) \mid \text{each $f_i(t)$ is monic, degree $d$, with no common root} \right\}. \]
Part of his proof included exhibiting homological stability via (non-algebraic) maps
\[  \Rat_d^n(\CC) \hookrightarrow  \Rat_{d+1}^n(\CC). \]
This work inspired many generalizations in the 40 years since, for example extending the domain to genus $g \geq 1$ curves and the target to Grassmannians or toric varieties; see, for example, \cite{BHM, CCMM, GKY, Guest, KM, Kirwan}. Jun-Yong Park remarked that his joint work with Changho-Han and Hunter Spink \cite{HP,PS} showed stability for maps from $\PP^1$ to weighted projective spaces $\PP(a,b)$.  Kirsten Wickelgren commented that the Chow ring of stable maps from $\PP^1 \to \PP^n$ is known to stabilize in the degree, due to Pandharipande \cite{pandharipande}. 

Knight asked what was known about extensions of Segal's results to the case
\[ \Hol_d(\CC\PP^m,\CC\PP^n) = \left\{ (f_0,\dots,f_n) : \text{each } f_i \in \CC[x_0,\dots,x_m]_d \text{ with no common root} \right\}, \]
where $1 < m \leq n$ and $d \geq 1$. Claudio G\'omez-Gonz\'ales remarked that this work had been carried out, first by Mostovoy \cite{Mostovoy03, Mostovoy12} and, more generally, by Mostovoy--Munguia-Villanueva \cite{MMV}. Further,  \cite{gomezgonzales} calculates the stable homology as $d \to \infty$, and shows that both the unstable rational homology and corresponding point counts seem to be independent of $d$ (except for dimension reasons, for the latter). Benson Farb added that, in general, little is known about spaces of algebraic maps when the (complex) dimension of the domain is greater than one (except when those spaces are finite).  This motivates: 

\begin{question}
	Can we develop more general methods to study the topology of spaces of algebraic maps $X \to Y$ between varieties, especially when $X$ has dimension greater than 1?
\end{question}

\begin{question}
	Can we replace homology/cohomology with the Chow ring in more general settings?
\end{question}

%20
\subsection{Jesse Kass}

\begin{question} 
What are the standing obstructions to computing cohomologies of discriminants? Are there problems that can serve as benchmarks for extending the existing techniques?
\end{question}

Orsola Tommasi commented that she uses Vassiliev's method (see \cite{vassiliev, tommasi}). Many papers cited here use some variant of Vassiliev's spectral sequence; Mostovoy \cite{Mostovoy12}, for example, use a truncated version of the spectral sequence when the limiting terms become too cumbersome to compute.

%21
\subsection{Jun-Yong Park}

Trace formulas give a concrete way to work with arithmetic topology. However, moduli functors that we want to study are often represented in the category of stacks rather than schemes, which are harder to work with. Kai Behrend \cite{behrend} established a trace formula for algebraic stacks, which has extended this arithmetic topology bridge to many moduli stacks. For example, Han and Park \cite{HP,HP2,park20} computed the $\ell$-adic cohomology of the Hom stack $\Hom_n(\PP^1, \PP(\lambda_0, \ldots, \lambda_N))$ and its class in the Grothendieck ring, with connections to fine moduli stacks of elliptic and hyperelliptic fibrations over $\PP^1$.

\begin{question}
	Can we develop a robust and general theory of point counting to understand moduli stacks arithmetically?
\end{question}

Wei Ho remarked on important connections to the work of Ellenberg--Satriano--Zureick-Brown \cite{ESZB}, bringing together many arithmetic problems about counting points, such as the Batyrev--Manin Conjecture, the Malle Conjecture, and more via, among other things, establishing a theory of heights on stacks. See, for example, the paper by Boggess and Sankar \cite{BS}, which discusses stacky versions of these conjectures in the context of the Ellenberg--Satriano--Zureick-Brown framework.

%22
\subsection{Inna Zakharevich}

Consider the Grothendieck ring of varieties over a field $k$, $K_0(\operatorname{Var}_k)$, equipped with its dimension filtration, where the $n$th graded piece of the filtration is given by those elements that can be represented as sums of varieties of dimension at most $n$, i.e., the image of the natural map 
\begin{equation}\label{natural map}
	K_0(\operatorname{Var}_k^{\leq n}) \coloneqq \ZZ[X| \dim X \leq n]/([X] = [Y]+[X-Y]) \to K_0(\operatorname{Var}_k). \end{equation}
The quotients $K_0(\operatorname{Var}_k^{\leq n})/K_0(\operatorname{Var}_k^{\leq n-1})$ have a simple presentation as the free abelian group on the birational isomorphism classes of $n$-dimensional varieties. 

However, as shown by Karzhemanov \cite{karzhemanov} in 2014 and Borisov \cite{borisov} in 2015, the map (\ref{natural map}) is non-injective when $k \subseteq \CC$. The failure of injectivity is measured by a spectral sequence whose differentials compute the differences \[[X-U]-[X-V],\] where $U$ and $V$ are isomorphic open subsets identified by a birational automorphism of $X$; for details, see \cite{zakharevich}. This motivates the following problem:

\begin{problem} Fix a field $k$, variety $X$, and a birational automorphism $f: X \dashrightarrow X$. Let $U, V \subseteq X$ be open subsets such that $f: U \to V$ is an isomorphism. Give an explicit example such that the complements $X-U$ and $X-V$, which will be equal in the Grothendieck ring, are not piecewise isomorphic. 
\end{problem}

As remarked above, such an example must exist for $k \subseteq \CC$, though none is known; it is unknown if examples exist over other fields.

%23
\subsection{Zinovy Reichstein}
\begin{problem} Let $G$ be a finite group and fix subgroups $H_1,\dots,H_r \leq G$, not necessarily distinct. Moreover, fix $d > 0$ and faithful $d$-dimensional representations $V_i$ of each of the $H_i$. Is there a $d$-dimensional irreducible variety $X$ over $\CC$ with a $G$ action such that there are $r$ distinct points $x_1, \dots, x_r$ with $x_i$ fixed by $H_i$ and isotropy representation given by $V_i$?
\end{problem}

Reichstein remarked that this is analogous in spirit to the Chinese Remainder Theorem and can be done if each $H_i$ is abelian; see \cite[Theorem 8.6]{RY} for details.  Moreover, when $r=1$ and $H=G$, one can simply take $X = V$, while for $r=2$ and $H_1=H_2=G$, it is an unpublished result of A. Kresch. In the special case, where $H_1,\ldots, H_r$ are the non-conjugate Sylow subgroups of $G$ (i.e., a list of Sylow subgroups for the different primes dividing $|G|$), this problem arose in the context of work on essential dimension of finite groups (see [Conjecture 11.5]\cite{DR}).

\subsection{Hunter Spink}

\begin{problem}\label{spink_fraction} Let $R$ be a ring with $m_1, \dots, m_k \in R$. Assume that, for any subset $A \subset \{1,\dots,k\}$, we have
\[ \sum_{i \in A} m_i \in R^\times. \]
Also set $a_i = m_1 x_1^i + \dots + m_k x_k^i \in R[x_1,\dots,x_k]$. What is the smallest $n$ such that there exist monic $f(z), g(z)$ with coefficients in $R[a_1,\dots,a_n]$ giving
\[ \frac{f(z)}{g(z)} = \frac{m_1}{z-x_1} + \cdots + \frac{m_k}{z-x_k}? \]
\end{problem}

For context on this problem, see the paper by Spink and Dennis Tseng \cite{spinktseng} on incidence strata of affine varieties that also outlines their conjecture of $n = 2^k+1$ to problem \ref{spink_fraction}.

\bibliographystyle{amsplain} 
\bibliography{problems}

\end{document}